\documentclass[12pt,reqno]{amsart}
\usepackage[a4paper,margin=28mm]{geometry}
\usepackage{amsmath,amssymb,amsthm,mathtools,booktabs}
\usepackage{placeins}
\usepackage{xcolor}
\usepackage[colorlinks=true,linkcolor=blue,citecolor=blue,urlcolor=blue]{hyperref}

\makeatletter
\renewcommand\subsection{\@startsection{subsection}{2}%
  \z@{.8\linespacing\@plus.3\linespacing}{.35\linespacing}%
  {\normalfont\bfseries}}
\makeatother
\allowdisplaybreaks[2]
\newcommand{\C}{\mathbb C}
\newcommand{\mm}{\mathfrak m}
\newcommand{\gr}{\operatorname{gr}}
\newcommand{\ord}{\operatorname{ord}}
\newcommand{\Hilb}{\operatorname{Hilb}}
\newcommand{\rank}{\operatorname{rank}}

\newcommand{\tbin}[2]{\binom{#1}{#2}_{\!+}}
\newcommand{\Pclass}{\mathcal P_{N;m,d}}
\newcommand{\Pall}{\mathcal P_{N;m}}
\newcommand{\OO}{\mathcal O}
\newtheorem{theorem}{Theorem}[section]
\newtheorem{proposition}[theorem]{Proposition}
\newtheorem{lemma}[theorem]{Lemma}
\newtheorem{corollary}[theorem]{Corollary}
\newtheorem{thmintro}{Theorem}

\theoremstyle{definition}
\newtheorem{example}[theorem]{Example}
\newtheorem{remark}[theorem]{Remark}
\newtheorem{notation}[theorem]{Notation}
\newtheorem{definition}[theorem]{Definition}
\title[Sharp minima of higher Milnor and Tjurina numbers]
{Sharp minima of higher Milnor and Tjurina numbers of isolated singularities}
\author{Chuangqiang Hu}
	\address{Sun Yat-Sen University, School of Mathematics, Guangzhou, 510275, P. R. China} 
	\email{huchq5@mail.sysu.edu.cn}
	\thanks{Hu was supported by NSFC Grant 12441107.}
	
	\author{Stephen S.-T.~Yau}
	\address{Beijing Institute of Mathematical Sciences and Applications,  Beijing 101408, P. R. China} 
	\email{yau@uic.edu}
	\thanks{Corresponding author: Stephen~S.-T. Yau\\
	Yau was supported by Tsinghua University Education Foundation fund (042202008).}
  
  \author{Lishan Yu}
  \address{Beijing Institute of Mathematical Sciences and Applications, Beijing, 101408, P.~R.~China}
  \email{lishany@bimsa.cn}
	
	\author{Huaiqing~Zuo}
	\address{Department of Mathematical Sciences, Tsinghua University, Beijing 100084, P. R. China}
	\email{hqzuo@mail.tsinghua.edu.cn}
	\thanks{Zuo was supported by NSFC Grant 12671056.}


\begin{document}
\begin{abstract}
We determine the minimum higher Milnor and Tjurina numbers of isolated
complex polynomial singularities with fixed ambient dimension and multiplicity.
A Fermat polynomial of multiplicity $m$, perturbed by the $(m+1)$-st power
of a linear form, attains both minimum sequences at every order.
These values remain the minima at each prescribed total degree $d>m$;
the homogeneous case $d=m$ has a separate Tjurina formula.
The proof expresses the Tjurina defect through multiplication in the
filtered Milnor algebra and uses the monomial strong Lefschetz theorem
to attain the resulting rank bound. We also determine the associated
graded higher Tjurina algebras of this family and show that they are
unchanged by terms of order at least $m+2$.
\end{abstract}
\maketitle
\tableofcontents

\section{Introduction and the minimum theorem}

\subsection{The invariants and the prescribed class}

The aim of this paper is to determine, for fixed ambient dimension and
multiplicity, the least possible values of all higher Milnor and
Tjurina numbers of polynomials with an isolated critical point at
the origin. One polynomial of degree $m+1$ attains both entire minimum sequences.
Adding higher-degree terms gives an attaining polynomial of each prescribed
degree $d>m$.

Put $\OO=\C[\![x_1,\ldots,x_N]\!]$. Let $\mm=(x_1,\ldots,x_N)$ be the maximal ideal and denote by 
\[ 
 J(f)=(f_{x_1},\ldots,f_{x_N})
\]
the Jacobian ideal of $f$. 
 \begin{definition}\label{def:inv}
For a polynomial $f$ with $f(0)=0$ and an isolated critical point
at the origin, define

\[
\mu_k(f)=\dim_\C \OO / \mm^k J(f),
\]
and
\[
\tau_k(f) = \dim_\C \OO/ (f,\mm^k J(f)), \qquad  k\geq 0.
\]
\end{definition}
  
At $k=0$ these are the usual Milnor and Tjurina numbers. For $k>0$,
replacing $J(f)$ by the smaller ideal $\mm^kJ(f)$ retains additional
classes in the quotients. We work in the completed local ring; completion
preserves these finite lengths, so the convergent local ring gives the
same invariants.

Fix integers $N\geq1$ and $m\geq2$. The class for the
main theorem is
\begin{equation}\label{eq:unrestrictedclass}
 \Pall=
 \{f\in\C[x_1,\ldots,x_N]:
       \ord_0f=m,\ \dim_\C \OO/J(f)<\infty\}.
\end{equation}
 The condition on $J(f)$ concerns only the germ at the origin.
The main theorem imposes no restriction on total degree. For the
exact-degree refinement, we use the subclass
$\Pclass=\{f\in\Pall:\deg f=d\}$, where $d\geq m$.

\subsection{The minimum formulas}

The minimum formulas involve two functions: the homogeneous Milnor value
$B_{N,m}(k)$ and a correction $Q_{N,m}(j)$ accounting for the further
decrease in the Tjurina number caused by nonhomogeneous terms.
We define both from the Hilbert function of a monomial complete intersection.

\begin{notation}\label{not:parameters}
Put $q=m-1$ and $S=N(m-2)$, and define $h_i$ by
\begin{equation}\label{eq:parameters}
 \sum_{i=0}^{S}h_it^i=(1+t+\cdots+t^{q-1})^N.
\end{equation}
We set $h_i=0$ outside $0\leq i\leq S$.
\end{notation}
Thus $h_i$ is the dimension in degree $i$ of the monomial complete
intersection
\[
 \C[x_1,\ldots,x_N]/(x_1^q,\ldots,x_N^q).
\]
\begin{proposition}\label{prop:parameters}
If $f$ is homogeneous of degree $m$ and has an isolated critical
point at the origin, then $\mu_k(f)=B_{N,m}(k)$ for every $k\geq0$,
where $B_{N,m}(k)$ is given by
\begin{equation}\label{eq:B}
 B_{N,m}(k)=\binom{N+m+k-2}{N}
                    +\sum_{i=m+k-1}^{S}h_i.
\end{equation}
\end{proposition}

The sum in \eqref{eq:B} gives the total dimension of the graded
components of the Milnor algebra in degrees at least $m+k-1$.
It vanishes when $m+k-1>S$.

This is the homogeneous case of the formula in
\cite[Proposition~4.1]{LiuYauZhang2026}.

At order zero, \eqref{eq:B} gives $B_{N,m}(0)=(m-1)^N$.

\begin{notation}
For the Tjurina correction, set

\begin{equation}\label{eq:Q}
 Q_{N,m}(j)=\sum_{i\geq j}\min\{h_i,h_{i+m+1}\},\qquad j\geq0.
\end{equation}
\end{notation}
The sum is finite: only indices $j\leq i\leq S-m-1$ can
contribute, and the correction vanishes for $j>S-m-1$. In
particular, $Q_{N,m}$ is identically zero when $S<m+1$.
Finite binomial expressions for both $B_{N,m}$ and $Q_{N,m}$ are
given in Section~\ref{sec:explicit}.

To cover the small values of $k$, we use the truncated binomial
convention
\[
 \tbin{u}{b}=
 \begin{cases}\binom{u}{b},&u\geq b,\\0,&u<b.\end{cases}
\]
Here $b$ is a nonnegative integer. In particular,
$\tbin{N+k-2}{N}$ vanishes for $k=0,1$.

\begin{thmintro}\label{thm:main}
For every $N\geq1$, $m\geq2$, and $k\geq0$:
\begin{enumerate}
\item The minimum Milnor numbers are
\begin{equation}\label{eq:minmu}
 \min_{f\in\Pall}\mu_k(f)=B_{N,m}(k).
\end{equation}
\item The minimum Tjurina numbers are
\begin{align}
 \min_{f\in\Pall}\tau_0(f)&=(m-1)^N-Q_{N,m}(0),
       \label{eq:minnu0}\\
 \min_{f\in\Pall}\tau_k(f)&=B_{N,m}(k)
              -\tbin{N+k-2}{N}-Q_{N,m}(k-1),\qquad k\geq1.
       \label{eq:minnuk}
\end{align}
% \item For each fixed $N,m$, the polynomial
% \[
%  F_{N;m}:=\sum_{i=1}^N x_i^m+(x_1+\cdots+x_N)^{m+1}
% \]
% attains all these minima simultaneously for every $k\geq0$.
\end{enumerate}
\end{thmintro}

\begin{thmintro}\label{thm:family}
Put $\ell=x_1+\cdots+x_N$. For every $\lambda\in\C^\times$ and
every polynomial $g\in\mm^{m+2}$, the polynomial
\begin{equation}\label{eq:family}
 f=\sum_{i=1}^N x_i^m+\lambda\ell^{m+1}+g,\qquad
 A=\OO/J(f),\qquad a=[f]
\end{equation}
belongs to $\Pall$ and attains all the minima of
Theorem~\ref{thm:main} simultaneously at every order $k\geq0$.
We call such an $f$ a \emph{Fermat--Lefschetz polynomial}.
\end{thmintro}

\begin{remark}\label{rem:seriestails}
Theorem~\ref{thm:family} also holds for a formal tail
$g\in\mm^{m+2}$. The proof in Section~\ref{sec:lefschetz} takes
place in $\OO$ and shows that the resulting series has an isolated
critical point at the origin and that its higher Milnor and Tjurina
numbers still equal the minima of Theorem~\ref{thm:main}. The
polynomial hypothesis is needed only for membership in $\Pall$.
\end{remark}

Taking $g=0$ gives a simultaneous attainer of degree $m+1$, and
taking $g=x_1^d$ gives one of each degree $d\geq m+2$.

\begin{thmintro}\label{thm:generic}
Let $V_{m,d}=\bigoplus_{v=m}^{d}\C[x_1,\ldots,x_N]_v$ with $d > m $.
There is a nonempty Zariski open subset of $V_{m,d}$ contained
in $\Pclass$ on which both entire sequences $(\mu_k)_{k\geq0}$
and $(\tau_k)_{k\geq0}$ attain the minimum values in
Theorem~\ref{thm:main}.
\end{thmintro}

\subsection{Relation to earlier work}

The algebra structures underlying Definition~\ref{def:inv} also enter
the classification of singularities. By the Mather--Yau theorem \cite{MatherYau}, the Tjurina algebra
determines the analytic isomorphism class of an isolated complex
hypersurface singularity. Greuel--Pham \cite{GreuelPham} obtain
a determination theorem in positive characteristic using higher
Tjurina algebras and also study higher Milnor algebras equipped with
the action of the defining function. Here we ask how small the lengths
of these higher algebras can be when only dimension and multiplicity
are prescribed.

 At order zero, the Milnor bound
$\mu_0(f)\geq(m-1)^N$ and its semihomogeneous equality case are
special cases of the weighted result of Furuya--Tomari
\cite{FuruyaTomari}; see also
Yau--Zuo \cite[Theorem~2.5 and Proposition~A]{YauZuo2012}.
Here semihomogeneous means that the degree-$m$ part has isolated
gradient. The role of homogeneity in the Tjurina number is also
classical: Saito \cite{Saito} characterizes the equality
$\mu_0=\tau_0$ by weighted homogeneity after an analytic change of
coordinates. Yau--Zuo \cite[Corollary~B]{YauZuo2015} characterize
analytic equivalence to a homogeneous singularity by
$\mu_0=\tau_0=(m-1)^N$. Thus the ordinary Milnor minimum and the
homogeneous ordinary Tjurina value in Remark~\ref{rem:homogeneous}
are established results. 

 For the higher sequences, Hu--Yau--Zuo
\cite[Corollary~4]{HYZHomogeneous} obtain the homogeneous Tjurina
formula in \eqref{eq:homogeneous}. Hussain--Liu--Yau--Zuo
\cite{HLYZ} study both sequences for weighted homogeneous
singularities and formulate a lower-bound problem at fixed weight
type. Hu--Yau--Zuo \cite{HYZ} develop a Koszul-complex method for
computing the higher Tjurina series and carry out the
three-variable weighted homogeneous calculation. Further work of
Hussain--Zuo \cite{HussainZuo2026} relates higher moduli algebras
and their derivation Lie algebras, including computations for simple
singularities at small orders.\par

Shi--Yau--Zuo \cite[Theorem~4.63 and Conjecture~4.65]{ShiYauZuo2024}
prove invariance of $\mu_k$ under semiquasihomogeneous perturbations for
$k\leq\min_i\ord(f_{x_i})$ and ask whether it holds at every order.
For ordinary semihomogeneous germs in $N\geq2$ variables,
Liu--Yau--Zhang \cite[Main Theorem~A and Proposition~4.1]{LiuYauZhang2026}
establish this all-order invariance and the formula $B_{N,m}(k)$.
Their Proposition~4.3 determines the optimal uniform range in which
the higher Tjurina numbers equal the homogeneous values.

Ribeiro--Fern\'andez-P\'erez \cite{FernandezPerezRibeiro} extend
the higher invariants to plane foliations. For plane function germs,
Corollary~4.1 of that paper expresses every $\mu_k$ in terms of
$\mu_0$ and the multiplicity; together with the ordinary Milnor bound,
it gives the $N=2$ case of our minimum Milnor sequence.
Their Theorem~6.1 gives a higher Tjurina lower bound for reduced
weighted homogeneous plane curves of normalized weight type
$(w_1,w_2;1)$ satisfying $\ord_0 f\geq1/w_1$.

Theorem~\ref{thm:main} covers arbitrary tangent cones and nonhomogeneous
terms in every dimension.

Among semiquasihomogeneous deformations with principal part
$x^a+y^b$, Brian\c{c}on, Granger, and Maisonobe \cite{BGM}
compute the minimum ordinary Tjurina number recursively; see
Almir\'on--Blanco \cite[\S3, Theorem~2]{AlmironBlanco}
for a restatement.

For ordinary plane multiple points, this yields the value
$\lfloor(3m^2-2m-4)/4\rfloor$. Canino--Gimigliano--Id\`a
\cite[\S5 of the preprint version]{CaninoGimiglianoIda}
give an elementary calculation and explicit degree-$(m+1)$
polynomials attaining it. We recover this classical value in
\eqref{eq:plane}.

Beyond ordinary multiple points, the ordinary Tjurina minimum is also
studied within a fixed topological type. For irreducible plane curves,
Alberich-Carrami\~nana and coauthors \cite{AlberichEtAl} express the
minimum in terms of resolution data, while Genzmer--Hernandes
\cite{GenzmerHernandes} obtain formulas through Saito bases.

For reducible plane curves, Genzmer
\cite{Genzmer2022,Genzmer2024} develops an algorithm for the generic
dimension of the equisingular moduli space, first for unions of
smooth branches and then for arbitrary reduced plane curve germs.
Almir\'on--Hernandes \cite[Theorem~4.1]{AlmironHernandes2026} treat two

equisingular branches with sufficiently large intersection
multiplicity and show that the Tjurina number is constant on the
corresponding equisingularity class.  
These equisingular problems keep the topological type fixed.
In the plane class $\mathcal P_{2;m}$, both the number of branches
and the topological type may vary.

The closest higher-dimensional ordinary calculation is
Wahl's minimum for positive-weight deformations of the
three-variable Fermat polynomial \cite[Example~4.7]{Wahl}.
Eriksson \cite[Theorem~2.2 and Proposition~3.2]{Eriksson}
gives a lower bound over all isolated three-variable germs of the
prescribed multiplicity and the one-term attainer
\[ 
x^m+y^m+z^m+(x+y+z)^{m+1}.\]
Formula~\eqref{eq:surface} recovers this three-variable minimum.
In arbitrary dimension,

Ma--Zuo \cite[Theorem~5.3 and its proof, Steps~II--III]{MaZuo}
compute the ordinary Tjurina Hilbert function of general
positive-weight Fermat deformations and deduce the corresponding
asymptotic minimum of the ordinary Tjurina number.

 Their finite-degree Hilbert-function calculation already
contains the ordinary correction used in \eqref{eq:minnu0}.

The Milnor part of Theorem~\ref{thm:main} follows from the
semihomogeneous formula by specialization.
For the Tjurina numbers, we bound the images of multiplication by
$[f]$ in the filtered Milnor algebra of an ordinary semihomogeneous germ.
This estimate yields the correction $Q_{N,m}$ and, by specialization,
lower bounds for every tangent cone. The family in
Theorem~\ref{thm:family} attains these bounds simultaneously.

The rank calculation also determines the associated

graded higher Tjurina algebras of the attaining family
(Corollary~\ref{cor:gradedattainers}), which are unchanged by tails
of order at least $m+2$. This algebraic conclusion is related to the
use of Tjurina Hilbert functions to distinguish deformation strata
in the work of Greuel--Hertling--Pfister \cite{GreuelHertlingPfister}.

\subsection{Structure and proof strategy}

Suppose first that the degree-$m$ part of $f$ has isolated gradient.
Let $A=\OO/J(f)$, with maximal ideal
$\mm_A$, and put $a=[f]\in A$.
Euler's identity removes the degree-$m$ term from this class, so
$a\in\mm_A^{m+1}$. This is the source of the shift by $m+1$ in
$Q_{N,m}$.
Section~\ref{sec:filtrations} expresses the extra Tjurina decrease
as $\dim_\C(a\mm_A^{k-1})$ for $k\geq1$; at $k=0$ it is
$\dim_\C(aA)$. Thus the minimum problem becomes an upper-bound
problem for multiplication on a filtered algebra.

The Hilbert function $(h_i)$ bounds these multiplication images.
For the family \eqref{eq:family}, the monomial strong Lefschetz theorem
\cite{Stanley,PhuongTran} gives the corresponding maximal ranks in
the associated graded algebra. Lifting independent products to the
local algebra gives the opposite inequality, so the filtered bound
is attained at every order. Sections~\ref{sec:filtrations}
and~\ref{sec:lefschetz} carry out these two steps.

Section~\ref{sec:fullclass} extends the bounds to arbitrary tangent cones
by deformation and finite-jet semicontinuity. For generic attainment
at a fixed degree $d>m$, a uniform nilpotence bound leaves only finitely
many rank conditions. The corresponding open loci contain a common
member of the explicit family. The required Tjurina cutoff on the ordinary semihomogeneous
locus is also established in
\cite[Proposition~4.3]{LiuYauZhang2026}.
Euler's identity gives the homogeneous boundary values.

\section{The Tjurina defect and a filtered rank bound}
\label{sec:filtrations}

\subsection{The Jacobian filtration}
\begin{definition}
Write $f=f_m+f_{m+1}+\cdots$, with each $f_v$ homogeneous of degree
$v$. We call $f$ \emph{ordinary semihomogeneous} if $f_m$ has
isolated gradient.
\end{definition}
For $N\geq2$, this means that the projective hypersurface $f_m=0$
is smooth in $\mathbb P^{N-1}$. The partial derivatives of $f_m$
all have the same degree
$q=m-1$ and form a regular sequence. This regularity will control
the ideals $\mm^kJ(f)$ even after higher terms are added to $f$.

For a quotient $A=\OO/I$ with maximal ideal $\mm_A$, the associated
graded algebra is
\[
 \gr A=\bigoplus_{v\geq0}\mm_A^v/\mm_A^{v+1}.
\]
Its dimension in degree $v$ is the corresponding Hilbert-function
value, and $\Hilb(\gr A,t)$ denotes the series of these dimensions.

The initial form $g^*$ of a nonzero $g\in \OO$ is its lowest nonzero
homogeneous part. We have $\gr \OO=P$, where
$P=\C[x_1,\ldots,x_N]$. The algebra $\gr(\OO/I)$ is obtained by
quotienting $P$ by the ideal generated by the initial forms of all
elements of $I$.
If the initial forms of the chosen generators form a regular
sequence, they generate the initial ideal of $I$.

The following equal-order form of the Valabrega--Valla criterion
\cite[Corollary~2.7]{ValabregaValla} relates the intersections
$I\cap\mm^{q+k}$ to the ideals $\mm^kI$.

For the Jacobian case, see also
\cite[Proposition~3.1]{LiuYauZhang2026}.

\begin{lemma}\label{lem:strict}
Let $g_1,\ldots,g_N\in \OO$ have initial forms $g_1^*,\ldots,g_N^*$
of degree $q$ forming a regular sequence in
$P=\C[x_1,\ldots,x_N]$. For $I=(g_1,\ldots,g_N)$:
\begin{enumerate}
\item $\gr_{\mm}(\OO/I)\simeq P/(g_1^*,\ldots,g_N^*)$;
\item \label{it:strict} $I\cap\mm^{q+k}=\mm^kI$ for every $k\geq0$.
\end{enumerate}
\end{lemma}
\begin{proof}
Apply \cite[Corollary~2.7]{ValabregaValla} to the $\mm$-adic
filtration and the generators of common order $q$. It gives
\[
 I\cap\mm^v=\sum_{i=1}^N\mm^{v-q}g_i\qquad(v\geq q).
\]
Taking $v=q+k$ proves~\ref{it:strict}. If $b\in I$ has order $v$,
this representation expresses its initial form as a combination
of the $g_i^*$. Thus the initial ideal of $I$ is
$(g_1^*,\ldots,g_N^*)$, proving~(i).
The criterion applies in the complete local ring $\OO$.
\end{proof}

We use the semihomogeneous Milnor formula
\cite[Propositions~3.1 and~4.1]{LiuYauZhang2026} in the following
form. The short length calculation also applies to formal series
and to $N=1$.

\begin{proposition}\label{prop:ordinarymu}
Let $f$ be ordinary semihomogeneous of multiplicity $m$, and put
$A=\OO/J(f)$, with maximal ideal $\mm_A$. Then:
\begin{enumerate}
\item \begin{equation}\label{eq:ordinaryA}
 \Hilb(\gr A,t)=(1+t+\cdots+t^{q-1})^N;
\end{equation}
\item $\dim_\C A=q^N$ and $\mm_A^{S+1}=0$;
\item $\mu_k(f)=B_{N,m}(k)$ for every $k\geq0$.
\end{enumerate}
\end{proposition}
\begin{proof}
\emph{(i), (ii).}
Lemma~\ref{lem:strict} gives $\gr A=P/J(f_m)$.
The $N$ initial partials form a regular sequence of degree $q$,
so this complete intersection has Hilbert series
$(1-t^q)^N/(1-t)^N=\sum_i h_it^i$.
The Hilbert series has no terms above degree $S$, so
\[
 \mm_A^{S+1}/\mm_A^{S+2}=0,\qquad
 \mm_A^{S+1}=\mm_A\mm_A^{S+1}.
\]
Nakayama's lemma gives $\mm_A^{S+1}=0$.
Summing the graded dimensions now gives
$\dim_\C A=\sum_i h_i=q^N$.

\emph{(iii).}
For the last claim, strictness identifies
$\mm^{q+k}/\mm^kJ(f)$ with $\mm_A^{q+k}$. Hence
\[
 0\longrightarrow\mm_A^{q+k}
 \longrightarrow\OO/\mm^kJ(f)
 \longrightarrow\OO/\mm^{q+k}\longrightarrow0
\]
is exact. Its outer terms have lengths
$\sum_{i=q+k}^S h_i$ and $\binom{N+q+k-1}{N}$, respectively.
Their sum is $B_{N,m}(k)$.
\end{proof}

\subsection{The exact Tjurina defect of semihomogeneous singularities}

Passing from $\mu_k$ to $\tau_k$ imposes the additional relation
$f=0$. At order zero, this amounts to quotienting the Milnor algebra
$A$ by the ideal generated by $a=[f]$.

 The ordinary defect also has the description
\[ 
\mu_0-\tau_0=\dim_\C \OO/(J(f):f).\]
Dimca--Sticlaru \cite[\S3]{DimcaSticlaru2026} use this colon ideal
to compare the ordinary plane defect with Fitting ideals of the
Tjurina ideal. Here we need the images of multiplication by $a$
on successive powers of the maximal ideal of $A$. 

At higher orders, the coefficients in a Jacobian representation
also matter. For homogeneous $f$, Euler's identity uses coefficients
of degree one and gives $f\in\mm J(f)$. This does not imply
$f\in\mm^kJ(f)$ for $k\geq2$, so the higher Tjurina
defect can be nonzero even in the homogeneous case.

For an ordinary semihomogeneous polynomial, there may also be a
nonzero class $a$ in $A$, of order at least $m+1$.
Using the colon identity in
\cite[equation~(4.11)]{LiuYauZhang2026}, we express the additional term
as $\dim_\C(a\mm_A^{k-1})$ for $k\geq1$.
The following proposition separates this term from the homogeneous
contribution.

\begin{proposition}\label{prop:defect}
Let $f$ be as in Proposition~\ref{prop:ordinarymu}, and let $a=[f]\in A=\OO/J(f)$.
Then:
\begin{enumerate}
\item $a\in\mm_A^{m+1}$ and
\begin{equation}\label{eq:defect0}
 \tau_0(f)=q^N-\dim_\C(aA);
\end{equation}
\item for $k\geq1$,
\begin{equation}\label{eq:defectk}
 \tau_k(f)=B_{N,m}(k)-\tbin{N+k-2}{N}
                   -\dim_\C(a\mm_A^{k-1}).
\end{equation}
\end{enumerate}
\end{proposition}
\begin{proof}
\emph{(i).}
The Euler operator $E=\sum_i x_i\partial_{x_i}$ satisfies
$Ef\in J(f)$, so

\begin{equation}\label{eq:eulerclass}
 a=[f-Ef/m]=\sum_{v>m}(1-v/m)[f_v]\in\mm_A^{m+1}.
\end{equation}
At $k=0$, the quotient $\OO/(f,J(f))=A/aA$ gives
\eqref{eq:defect0}.

\emph{(ii).}
For $k\geq1$, consider multiplication by $f$ as a map from $\OO$ to
$\OO/\mm^kJ(f)$. Its image has dimension $\mu_k(f)-\tau_k(f)$,
and its kernel is the colon ideal $(\mm^kJ(f):f)$. Strictness gives

\begin{equation}\label{eq:colon}
 (\mm^kJ(f):f)=(J(f):f)\cap\mm^{k-1}.
\end{equation}
For $b\ne0$, strictness gives
$bf\in\mm^kJ(f)$ precisely when $bf\in J(f)$ and
$\ord(bf)\geq q+k$. Since $\ord f=m=q+1$, the latter
condition is $b\in\mm^{k-1}$.

Factoring out this kernel gives the exact sequence
\[
 0\longrightarrow \OO/(\mm^kJ(f):f)
 \xrightarrow{\ \cdot f\ }\OO/\mm^kJ(f)
 \longrightarrow \OO/(f,\mm^kJ(f))\longrightarrow0.
\]
Thus $\mu_k(f)-\tau_k(f)$ is the length of the left-hand quotient.
By \eqref{eq:colon}, its defining ideal is
$C_k=(J(f):f)\cap\mm^{k-1}$.

To separate the low-order coefficient contribution, use
$C_k\subseteq\mm^{k-1}$ to form a second exact sequence:
\[
 0\longrightarrow\mm^{k-1}/C_k
   \longrightarrow \OO/C_k
   \longrightarrow \OO/\mm^{k-1}\longrightarrow0.
\]
The monomials of degree below $k-1$ form a basis of the right-hand
quotient. Its length is $\tbin{N+k-2}{N}$, depending
only on $N$ and $k$.

The left-hand term records the contribution that depends on $f$.
To identify it inside $A$, consider the map
\[
 \mm^{k-1}\longrightarrow A,\qquad b\longmapsto[bf],
\]
which has kernel $C_k$ and image $a\mm_A^{k-1}$.
This identifies $\mm^{k-1}/C_k$ with $a\mm_A^{k-1}$. Taking
lengths gives
\[
 \mu_k(f)-\tau_k(f)
 =\tbin{N+k-2}{N}+\dim_\C(a\mm_A^{k-1}).
\]
Substituting $\mu_k(f)=B_{N,m}(k)$ proves \eqref{eq:defectk}.
\end{proof}

At $k=1$, the binomial term is zero, so the defect is $\dim_\C(aA)$,
just as at order zero. For larger $k$, the binomial term increases
while the dimension of $a\mm_A^{k-1}$ cannot increase. The remaining
task is to bound this multiplication image.

\subsection{A rank bound from the Hilbert function}

For ordinary semihomogeneous germs, Proposition~\ref{prop:defect}
reduces the Tjurina minimum problem to maximizing
$\dim_\C(a\mm_A^j)$. Here $a\in\mm_A^{m+1}$ and $A$ has Hilbert
function $(h_i)$.

Each summand $\min\{h_i,h_{i+m+1}\}$ in $Q_{N,m}(j)$ is the
largest possible rank of a linear map between graded components of
dimensions $h_i$ and $h_{i+m+1}$. The lower limit $j$ determines which
source degrees contribute. To use this sum, we must show that the
degreewise bounds also control the full image in the filtered
algebra.

The symmetry and unimodality of $(h_i)$ determine which of $h_i$
and $h_{i+m+1}$ is smaller: in the early degrees the smaller one is
$h_i$, and in the later ones it is $h_{i+m+1}$.

\begin{lemma}\label{lem:unimodal}
The sequence $(h_i)_{i=0}^S$ of Notation~\ref{not:parameters} is symmetric
and unimodal. Consequently, for $i\geq0$,
\[
 h_i\leq h_{i+m+1}\quad\hbox{if }2i+m+1\leq S,
\]
and
\[
 h_i\geq h_{i+m+1}\quad\hbox{if }2i+m+1\geq S.
\]

\end{lemma}
\begin{proof}
Symmetry follows by reversing each factor
$1+t+\cdots+t^{q-1}$.

For unimodality, induct on the number of factors. A single factor
has all coefficients equal to one. When $q=1$, the entire product
is the constant polynomial one. Suppose $(b_i)$ is the
symmetric unimodal sequence for $N-1$ factors, with support
$[0,S']$. Convolution with the next factor gives
\[
 h_i=b_i+b_{i-1}+\cdots+b_{i-q+1},\qquad
 h_i-h_{i-1}=b_i-b_{i-q}.
\]
If $i\leq(S'+q-1)/2$, the integer $i$ is no farther from the
center $S'/2$ than $i-q$, so $b_i\geq b_{i-q}$, with zero
extension outside the support. Thus $h_i$ is nondecreasing up
to the midpoint; symmetry gives the other half.

The stated inequalities follow by comparing distances from $S/2$.
\end{proof}

A product with zero class in degree $i+m+1$ may still be nonzero
in a higher filtration degree, so the graded multiplication maps
alone do not bound the full image. We split the filtration of
$\mm_A^j$ at one degree: below it, the image is bounded by the
dimension of its domain; from that degree onward, by the dimension
of the power of $\mm_A$ containing it.

\begin{lemma}\label{lem:rank}
Let $A$ be an Artinian local $\C$-algebra with maximal ideal
$\mm_A$ and Hilbert function $(h_i)$ from \eqref{eq:parameters}.
If $a\in \mm_A^{m+1} $, then
\begin{equation}\label{eq:rankbound}
 \dim_\C(a\mm_A^j)\leq Q_{N,m}(j)\qquad(j\geq0).
\end{equation}
\end{lemma}
\begin{proof}
Fix $j$ and an integer $c\geq j$. Choose a decomposition of
vector spaces
\[
 \mm_A^j=V\oplus\mm_A^c.
\]
The dimension of $aV$ is at most that of $V$, namely
$\sum_{i=j}^{c-1}h_i$. For the other summand, the order of $a$ gives
$a\mm_A^c\subseteq\mm_A^{c+m+1}$.
The decomposition of $\mm_A^j$ gives
\[
 a\mm_A^j=aV+a\mm_A^c.
\]
The two images may overlap, but the dimension of their sum is at
most the sum of the bounds just obtained. Thus
\begin{equation}\label{eq:splitrank}
 \dim_\C(a\mm_A^j)
 \leq\sum_{i=j}^{c-1}h_i+\sum_{i=c+m+1}^{S}h_i.
\end{equation}
To use the smaller Hilbert-function value in each range, choose
the splitting degree

\[
 c=\max\left\{j,\left\lceil\frac{S-m-1}{2}\right\rceil\right\}.
\]
For $i\geq j$, Lemma~\ref{lem:unimodal} identifies
$\min\{h_i,h_{i+m+1}\}$ with $h_i$ below $c$ and with $h_{i+m+1}$
from $c$ onward. Reindexing the second sum in
\eqref{eq:splitrank} therefore gives
\[
 \sum_{i=j}^{c-1}h_i+\sum_{v=c+m+1}^{S}h_v
 =\sum_{i=j}^{c-1}h_i+\sum_{i\geq c}h_{i+m+1}
 =\sum_{i\geq j}\min\{h_i,h_{i+m+1}\}=Q_{N,m}(j).
\]
If $m+1>S$, both the correction and the second sum vanish, and
$c=j$ gives the same conclusion.
\end{proof}

The argument uses only the Hilbert function and the order of $a$.
It therefore applies to every Milnor algebra in
Proposition~\ref{prop:ordinarymu}, regardless of its multiplication
table.

\begin{example}\label{ex:rank}
Take $N=3$ and $m=4$. Then $S=6$ and
$a\in\mm_A^5$.
Writing $A$ as the direct sum of the constants and $\mm_A$ gives
\[
 \dim_\C(aA)\leq1+\dim_\C\mm_A^6=2.
\]
The first term bounds the dimension of the image of the constants.
Every product of $a$ with an element of $\mm_A$ lies in
$\mm_A^6$, whose dimension gives the second term.
Restricting the source further gives
$\dim_\C(a\mm_A)\leq\dim_\C\mm_A^6=1$, whereas
$a\mm_A^2\subseteq\mm_A^7=0$.
These are the three bounds $Q_{3,4}(0)=2$, $Q_{3,4}(1)=1$, and $Q_{3,4}(2)=0$.
\end{example}

Combining the defect formula with Lemma~\ref{lem:rank} gives lower
bounds for the Tjurina numbers of ordinary semihomogeneous germs.
Section~\ref{sec:lefschetz} constructs polynomials attaining these bounds
at every order.

 \begin{proposition}\label{prop:semihomogeneous-values-bounds}
Let $f\in\Pall$ be ordinary semihomogeneous. Put
\[
A=\OO/J(f),\qquad a=[f]\in A.
\]
Then:
\begin{enumerate}
\item $\tau_0(f)=q^N-\dim_\C(aA)\ge q^N-Q_{N,m}(0)$;
\item for every $k\geq1$,
\begin{align*}
\tau_k(f)
& =
B_{N,m}(k)
-\tbin{N+k-2}{N}
-\dim_\C\big(a\mm_A^{k-1}\big) \\
& \ge
B_{N,m}(k)
-\tbin{N+k-2}{N}
-Q_{N,m}(k-1).
\end{align*}
\end{enumerate}
\end{proposition}
\begin{proof} 
By Proposition~\ref{prop:defect}, the two defects are expressed through
$aA$ and $a\mm_A^{k-1}$, with $a\in\mm_A^{m+1}$.
Lemma~\ref{lem:rank} bounds their dimensions by $Q_{N,m}(0)$ and
$Q_{N,m}(k-1)$, respectively. Substituting these bounds in
\eqref{eq:defect0} and \eqref{eq:defectk} proves the result.
\end{proof}

\section{Attainment by the Fermat--Lefschetz family}
\label{sec:lefschetz}

\subsection{Monomial Lefschetz ranks}

The graded Milnor algebra of the Fermat polynomial is
\[
 G=\C[x_1,\ldots,x_N]/(x_1^q,\ldots,x_N^q).
\]

The classical monomial strong Lefschetz theorem \cite{Stanley,PhuongTran}
gives the ranks of multiplication by powers of $\ell=x_1+\cdots+x_N$.
We recall its proof using an $\mathfrak{sl}_2$ action. 

\begin{lemma}[Strong Lefschetz theorem]\label{lem:slp}
 For every $s\geq0$ and $i\geq0$,
\begin{equation}\label{eq:slprank}
 \rank(\ell^s:G_i\longrightarrow G_{i+s})
          =\min\{\dim_\C G_i,\dim_\C G_{i+s}\},
\end{equation}
where $ G_i $ denotes the $i$-th component of $G$. 
\end{lemma}
\begin{proof}
We decompose $G$ into chains of consecutive degrees, all centered
at the middle degree. Multiplication by $\ell$ advances one step
along each chain. Their common center then makes it possible to
read injectivity and surjectivity from the endpoints.

An $\mathfrak{sl}_2$ action supplies the chains. In one variable,
define operators on $\C[x]/(x^q)$ by
\[
 \mathsf E(x^i)=x^{i+1},\qquad
 \mathsf F(x^i)=i(q-i)x^{i-1},\qquad
 \mathsf H(x^i)=(2i-q+1)x^i,
\]
where powers with exponents outside $0,\ldots,q-1$ are interpreted
as zero. These operators satisfy
\[
 [\mathsf H,\mathsf E]=2\mathsf E,\qquad
 [\mathsf H,\mathsf F]=-2\mathsf F,\qquad
 [\mathsf E,\mathsf F]=\mathsf H.
\]
On the tensor product of $N$ such spaces, the sums of these
operators give a finite-dimensional $\mathfrak{sl}_2$-module.
Here $\mathsf E$ is multiplication by $\ell$, and $\mathsf H$ acts on $G_i$
with weight $2i-S$.

By complete reducibility, $G$ is a direct sum of irreducible
finite-dimensional $\mathfrak{sl}_2$-modules over $\C$.
Each summand has one-dimensional weight spaces with weights
$-L,-L+2,\ldots,L$. Translating weight $2i-S$ into
degree $i$ places these spaces in the consecutive degrees
\[
 p,p+1,\ldots,S-p,\qquad p=(S-L)/2.
\]
If $v$ spans the lowest-degree space of one summand, that summand
has a basis
\[
 v,\mathsf Ev,\ldots,\mathsf E^{S-2p}v.
\]
Multiplication by $\ell$ moves along this chain by one degree.
All chains are centered at $S/2$, although they may begin at
different degrees $p$.

Suppose first that $2i+s\leq S$. A chain containing degree $i$
starts at some $p\leq i$ and ends at
\[
 S-p\geq S-i\geq i+s.
\]
It therefore also contains degree $i+s$, and $\mathsf E^s$ is
nonzero on its degree-$i$ space. Since this holds for every source
chain, the map $G_i\to G_{i+s}$ is injective.

Suppose instead that $2i+s\geq S$. A chain containing degree
$i+s$ ends at $S-p\geq i+s$, so
\[
 p\leq S-i-s\leq i.
\]
Thus this chain also contains degree $i$, and its degree-$(i+s)$
space is reached by $\mathsf E^s$. Every target chain is reached,
so the map is surjective. In either case its rank is the smaller
of the source and target dimensions, proving \eqref{eq:slprank}.
\end{proof}

\begin{example}\label{ex:slp}
The algebra $\C[x,y]/(x^3,y^3)$ has Hilbert function
$(1,2,3,2,1)$. For $\ell=x+y$, its Lefschetz chains have the
following shapes; a dot denotes a one-dimensional component:

\[
 \begin{array}{c|ccccc}
 \text{degree}&0&1&2&3&4\\ \hline
 \text{first chain}&\bullet&\bullet&\bullet&\bullet&\bullet\\
 \text{second chain}&&\bullet&\bullet&\bullet&\\
 \text{third chain}&&&\bullet&&
 \end{array}
\]
Multiplication by $\ell^2$ from degree one to degree three reaches
the target in both chains that occur in the source, so it has rank
two. From degree two to degree four, only the first chain reaches
the target; the rank is one, which is the full target dimension.
\end{example}

\subsection{Lifting the ranks to the local algebra}

\begin{proposition}\label{prop:attainment}
Let $\lambda\in\C^\times$, let $g\in\mm^{m+2}$, and let
$f$, $A$, and $a$ be as in \eqref{eq:family}.
Then:
\begin{enumerate}
\item $f$ has an isolated critical point at the origin;
\item $a$ is congruent to $-\frac{\lambda}{m}\,[\ell^{m+1}]$ modulo $\mm_A^{m+2}$; 
\item for every $j\geq0$,
\begin{equation}\label{eq:attainedrank}
 \dim_\C(a\mm_A^j)=Q_{N,m}(j).
\end{equation}
\end{enumerate}
\end{proposition}
\begin{proof}
\emph{(i).}
The initial partials are $mx_i^{m-1}$, so
Proposition~\ref{prop:ordinarymu} identifies $\gr A$ with $G$
in Lemma~\ref{lem:slp}; in particular, $A$ has finite length,
so the critical point at the origin is isolated.

\emph{(ii).}
For the polynomial \eqref{eq:family}, the Euler expression
\eqref{eq:eulerclass} cancels the degree-$m$ part and multiplies
$\lambda\ell^{m+1}$ by $-1/m$. The terms arising from $g$ have
order at least $m+2$, so

\begin{equation}\label{eq:aclass}
 a\equiv-\frac{\lambda}{m}\,[\ell^{m+1}]
             \pmod{\mm_A^{m+2}}.
\end{equation}

\emph{(iii).}
Fix $j\geq0$. We produce $\sum_{i\geq j}\rho_i$ independent vectors
in $a\mm_A^j$, where
\[
 \rho_i=\dim_\C\bigl(\ell^{m+1}G_i\bigr)=\min\{h_i,h_{i+m+1}\};
\]
the second equality is Lemma~\ref{lem:slp} with $s=m+1$.
Because $h_{i+m+1}=0$ for $i>S-m-1$, only the source degrees
$j\leq i\leq S-m-1$ contribute; for $j>S-m-1$ the sum is empty
and $Q_{N,m}(j)=0$.

For each contributing $i$, choose
$b_{i,1},\ldots,b_{i,\rho_i}\in G_i$ whose products with
$\ell^{m+1}$ form a basis of
$\ell^{m+1}G_i\subseteq G_{i+m+1}$.

By part (i), $\gr A=G$, so
$\mm_A^v/\mm_A^{v+1}\simeq G_v$ for every $v\geq0$.
Choose a lift $\widetilde b_{i,u}\in\mm_A^i$ of each $b_{i,u}$:
its class modulo $\mm_A^{i+1}$ corresponds to $b_{i,u}$ under this
isomorphism.

Since $a\in\mm_A^{m+1}$, each product $a\widetilde b_{i,u}$ lies in
$\mm_A^{i+m+1}$. Writing \eqref{eq:aclass} as
$a=-\lambda[\ell^{m+1}]/m+r$ with $r\in\mm_A^{m+2}$ gives
\[
 a\widetilde b_{i,u}
 =-\frac{\lambda}{m}\,\ell^{m+1}\widetilde b_{i,u}
   +r\widetilde b_{i,u},
 \qquad r\widetilde b_{i,u}\in\mm_A^{i+m+2}.
\]
The isomorphism $\gr A=G$ respects multiplication, so the class of
$a\widetilde b_{i,u}$ modulo $\mm_A^{i+m+2}$ corresponds to
$-\lambda\ell^{m+1}b_{i,u}/m$ in $G_{i+m+1}$.

We now prove that the products $a\widetilde b_{i,u}$
are independent. Suppose that
\[
 \sum_{i\geq j}\sum_{u=1}^{\rho_i}
       c_{i,u}\,a\widetilde b_{i,u}=0
\]
were a nontrivial relation, and let $i_0$ be the least source
degree having some nonzero coefficient $c_{i_0,u}$.
Every term with $i>i_0$ lies in
$a\mm_A^{i_0+1}\subseteq\mm_A^{i_0+m+2}$, so reduction modulo
$\mm_A^{i_0+m+2}$ leaves only the $i_0$ terms:
\[
 -\frac{\lambda}{m}\sum_{u}c_{i_0,u}\,\ell^{m+1}b_{i_0,u}=0
 \qquad\text{in }G_{i_0+m+1}.
\]
The products $\ell^{m+1}b_{i_0,u}$ are independent and
$\lambda\ne0$, so all $c_{i_0,u}$ vanish, contradicting the choice
of $i_0$. Hence no nontrivial relation exists.

All $\sum_{i\geq j}\rho_i$ products lie in $a\mm_A^j$, because
$i\geq j$ gives
$\widetilde b_{i,u}\in\mm_A^i\subseteq\mm_A^j$.
Their independence gives
\[
 \dim_\C(a\mm_A^j)\geq\sum_{i\geq j}\rho_i
 =\sum_{i\geq j}\min\{h_i,h_{i+m+1}\}=Q_{N,m}(j)
\]
by the definition \eqref{eq:Q}. Lemma~\ref{lem:rank} gives the
reverse inequality, so equality holds, proving
\eqref{eq:attainedrank}.
If $m+1>S$, then
$a\in\mm_A^{m+1}\subseteq\mm_A^{S+1}=0$ by
Proposition~\ref{prop:ordinarymu}, while $Q_{N,m}$ vanishes
identically as noted after \eqref{eq:Q}; both sides of
\eqref{eq:attainedrank} are then zero.

Since $j$ was arbitrary and the germ $f$ was fixed,
\eqref{eq:attainedrank} holds simultaneously for every $j\geq0$.
\end{proof}

\begin{remark}\label{rem:attainment-input}
The identification $\gr A=G$ and the congruence \eqref{eq:aclass}
hold for formal tails $g\in\mm^{m+2}$ as well as polynomial ones.
The same proof therefore establishes Proposition~\ref{prop:attainment}
for formal series.

\end{remark}

\begin{corollary}\label{cor:semihomogeneous-sharp}
The bounds in Proposition~\ref{prop:semihomogeneous-values-bounds}
are sharp. Every member of the Fermat--Lefschetz family
\eqref{eq:family} attains equality in all of them
simultaneously for every $k\ge 0$.
\end{corollary}

\begin{proof}
For each polynomial in the Fermat--Lefschetz family,
Proposition~\ref{prop:attainment} gives the multiplication ranks
at every order. Substituting these ranks into \eqref{eq:defect0}
and \eqref{eq:defectk} gives equality in both Tjurina bounds.

\end{proof}

%%%%%%%%%%%%%%%%%%%%%%%%%%%%%%%%%%%%%%%
% graded_algebras
%%%%%%%%%%%%%%%%%%%%%%%%%%%%%%%%%%%%%%%
\subsection{The associated graded Tjurina algebras}
\label{sec:graded}

For the Fermat--Lefschetz family, we can also describe
the associated graded Tjurina algebras by generators and relations.
These presentations determine the multiplication as well as the
Hilbert function.

At order zero, the Hilbert function below agrees with
Ma--Zuo \cite[proof of Theorem~5.3, Steps~II--III]{MaZuo}.

\begin{corollary}\label{cor:gradedattainers}
For the Fermat--Lefschetz family \eqref{eq:family}, put
$P=\C[x_1,\ldots,x_N]$, $\mathfrak n=(x_1,\ldots,x_N)\subset P$,
$f_m=\sum_i x_i^m$, and $J_0=(x_1^q,\ldots,x_N^q)$.
Then:
\begin{enumerate}
\item \begin{equation}\label{eq:gradedTzero}
 \gr_{\mm}\bigl(\OO/(f,J(f))\bigr)\simeq P/(J_0,\ell^{m+1});
\end{equation}
\item for $k\geq1$,
\begin{equation}\label{eq:gradedThigher}
 \gr_{\mm}\bigl(\OO/(f,\mm^kJ(f))\bigr)
   \simeq P/\bigl(f_m,\mathfrak n^kJ_0,
                         \ell^{m+1}\mathfrak n^{k-1}\bigr).
\end{equation}
\end{enumerate}
In particular, these graded algebras are independent of $\lambda\ne0$
and of every tail $g\in\mm^{m+2}$.
The Hilbert function in \eqref{eq:gradedTzero} is
$\max\{h_v-h_{v-m-1},0\}$ in degree $v$.
\end{corollary}
\begin{proof}
In each case, the proposed initial relations give a surjection from
the displayed homogeneous quotient onto the associated graded Tjurina
algebra. We prove that no further relations occur by comparing the two
finite lengths, using Propositions~\ref{prop:defect}
and~\ref{prop:attainment}.

\emph{(i).}
The initial partials show that the initial ideal of $(f,J(f))$
contains $J_0$. It also contains $\ell^{m+1}$, because
\[
 f-Ef/m=-\lambda\ell^{m+1}/m+
                  \text{terms of order at least }m+2.
\]
These containments induce a surjective graded map
\[
 P/(J_0,\ell^{m+1})\longrightarrow\gr\bigl(\OO/(f,J(f))\bigr).
\]
By Lemma~\ref{lem:slp}, the source $P/(J_0,\ell^{m+1})$ has
dimension $q^N-Q_{N,m}(0)$. Propositions~\ref{prop:defect}
and~\ref{prop:attainment} give the same length for the target.
The surjection is therefore between vector spaces of equal finite
dimension and has zero kernel, proving \eqref{eq:gradedTzero}.
The degree-$v$ component of this quotient has dimension
\begin{align*}
 h_v-\rank(\ell^{m+1}:G_{v-m-1}\longrightarrow G_v)
 & =h_v-\min\{h_{v-m-1},h_v\} \\
 & =\max\{h_v-h_{v-m-1},0\},
\end{align*}
where components of negative degree are zero.

\emph{(ii).}
Fix $k\geq1$. The initial ideal of $(f,\mm^kJ(f))$ contains
$f_m$ and $\mathfrak n^kJ_0$, coming from $f$ and the multiples
of its partials.

 To obtain the third family of relations, take a nonzero
homogeneous element $u\in\mathfrak n^{k-1}$. Then
\[
 u(f-Ef/m)=uf-\frac1m\sum_i ux_i f_{x_i}
                         \in(f,\mm^kJ(f)).
\]
The initial form of this element is a nonzero scalar multiple of
$u\ell^{m+1}$. Hence the initial ideal also contains
$\ell^{m+1}\mathfrak n^{k-1}$. We obtain a surjection
\[
 P/(f_m,\mathfrak n^kJ_0,\ell^{m+1}\mathfrak n^{k-1})
 \longrightarrow\gr\bigl(\OO/(f,\mm^kJ(f))\bigr).
\]
To compute the length of the source, first impose only the
homogeneous relations $f_m$ and $\mathfrak n^kJ_0$, and write
\[
 D_k=P/(f_m,\mathfrak n^kJ_0).
\]
The remaining relations form the image $U_k\subset D_k$ of the
ideal $\ell^{m+1}\mathfrak n^{k-1}$. Thus the proposed graded algebra
is $D_k/U_k$, and its length is obtained by subtracting the dimension
of this image from that of $D_k$.

For $D_k$, apply
Propositions~\ref{prop:ordinarymu} and~\ref{prop:defect}
to the homogeneous polynomial $f_m$. This gives
\[
 \dim_\C D_k =B_{N,m}(k)-\tbin{N+k-2}{N}.
\]
Euler's identity gives $f_m\in J_0$, so the multiplication term
in the defect formula vanishes.
To find $\dim_\C U_k$, we can work in the simpler algebra $G$.

The natural surjection
\[
 D_k\longrightarrow G=P/J_0
\]
is an isomorphism in each degree $v\geq q+k$.

In this range, the defining ideals of $D_k$ and $G$ have the same
homogeneous components. Indeed, each generator of $J_0$ has degree
$q$, so a degree-$v$ multiple has a coefficient of degree
$v-q\geq k$. Hence $(\mathfrak n^kJ_0)_v=(J_0)_v$.
Since $f_m\in J_0$, imposing it adds no relations in these degrees.

The ideal $\ell^{m+1}\mathfrak n^{k-1}$ starts in degree
$m+k=q+k+1$, so all components of $U_k$ lie in the range
where $D_k\to G$ is an isomorphism.

For each $i\geq k-1$, the degree-$(i+m+1)$ component of $U_k$ is
the image of $\ell^{m+1}P_i$ modulo $(J_0)_{i+m+1}$. Since
$\ell^{m+1}(J_0)_i\subseteq(J_0)_{i+m+1}$, this image is exactly
$\ell^{m+1}G_i$. In particular, the relevant map is
\[
 G_i\xrightarrow{\ \cdot\ell^{m+1}\ }G_{i+m+1},
\]
even when the source degree $i$ lies below the range in which
$D_k\to G$ is an isomorphism. Elements of $(J_0)_i$ contribute
zero to the target. Summing the ranks of the Lefschetz maps gives

\[
 \dim_\C U_k
 =\sum_{i\geq k-1}\rank(\ell^{m+1}:G_i\to G_{i+m+1})=Q_{N,m}(k-1).
\]
Consequently,
\[
 \dim_\C(D_k/U_k)
 =B_{N,m}(k)-\tbin{N+k-2}{N}-Q_{N,m}(k-1).
\]
By Propositions~\ref{prop:defect} and~\ref{prop:attainment}, this
is also the length of $\OO/(f,\mm^kJ(f))$. Equal finite dimensions
again force the surjection to have zero kernel, proving
\eqref{eq:gradedThigher}.
\end{proof}

\section{Proof of the minimum theorem}
\label{sec:fullclass}

For an arbitrary $f\in\Pall$, we deform its degree-$m$ part so that
the resulting homogeneous form has isolated gradient. General members
of this deformation have the same total degree as $f$ and are ordinary
semihomogeneous. Semicontinuity then transfers the lower bounds of
Section~\ref{sec:filtrations} to the original polynomial.

\subsection{Passing to a smooth tangent cone}

\begin{lemma}\label{lem:usc}
Given $f\in\Pclass$, there is a family $f_s$ of polynomials
such that
 \begin{enumerate}
  \item $f_0=f$;
  \item its general members are ordinary semihomogeneous and belong
to $\Pclass$;
\item for each fixed
$k\geq0$, the inequalities
\begin{equation}\label{eq:usc}
 \mu_k(f)\geq\mu_k(f_s),\qquad \tau_k(f)\geq\tau_k(f_s)
\end{equation}
hold on a nonempty Zariski open set of parameters $s$.
\end{enumerate}
\end{lemma}
\begin{proof}
 
\emph{(i), (ii).}
Let $g_m=\sum_i x_i^m$, and let $f_m$ be the degree-$m$ part of $f$.
Define the family
\[
 f_s=f+s(g_m-f_m).
\]
The degree-$m$ forms with isolated gradient form a nonempty Zariski
open set. To see this, consider the closed incidence set of pairs
$(h,[x])$ for which all partial derivatives of $h$ vanish at
$[x]\in\mathbb P^{N-1}$. Projection along the projective factor
is closed, so its image is the closed locus of forms without
isolated gradient. The complement contains $g_m$, whose partials
have no common projective zero.

At $s=1$, the degree-$m$ part of the pencil is $g_m$. The pencil
therefore meets this open set, so its degree-$m$ part has isolated
gradient for general parameters. For $N=1$, the condition reduces to
nonvanishing of the coefficient of $x_1^m$.

If $d>m$, the degree-$d$ part is unchanged; if $d=m$, general
members remain nonzero homogeneous forms of degree $m$.
The multiplicity and exact degree are consequently preserved.

\emph{(iii).}
Fix $k$. To compare lengths by matrix ranks, we need one finite
quotient that computes the lengths both at the original polynomial
and at the ordinary members of the family. Since $J(f)$ is
$\mm$-primary, choose $T$ with
$\mm^T\subseteq J(f)$ and put $L_0=T+k$; then
$\mm^{L_0}\subseteq\mm^kJ(f)$.

For the ordinary members, Proposition~\ref{prop:ordinarymu} gives
$\mm_{\OO/J(f_s)}^{S+1}=0$. Hence $\mm^{S+1}\subseteq J(f_s)$,
and multiplying by $\mm^k$ gives the
uniform inclusion
\[
 \mm^{S+1+k}\subseteq\mm^kJ(f_s).
\]
Choose $L\geq\max\{L_0,S+1+k\}$.
At $s=0$ and at every ordinary parameter, reduction modulo $\mm^L$
therefore preserves both lengths under consideration.

Use the monomials of degree below $L$ as a basis of $\OO/\mm^L$.
In this basis, form a matrix $M_\mu(s)$ whose columns are the
coordinate vectors of the monomial multiples
\[
 [x^\alpha(f_s)_{x_i}]\quad
 (1\leq i\leq N,\ k\leq|\alpha|<L-q).
\]
These columns span the image of $\mm^kJ(f_s)$. Adjoining the
columns $[x^\beta f_s]$ for $|\beta|<L-m$ gives $M_\tau(s)$,
whose full set of columns spans the image of $(f_s,\mm^kJ(f_s))$.
Multiples with $|\alpha|\geq L-q$ or $|\beta|\geq L-m$ vanish
in the chosen quotient: their orders are at least $L$ because
$\ord (f_s)_{x_i}\geq q$ and $\ord f_s\geq m$.

The same monomial basis and column indices are used for every $s$.
At $s=0$ and on the ordinary locus, the common truncation computes
the quotient lengths, giving the identities below. 
\[
 \begin{aligned}
 \mu_k(f_s)&=\dim_\C(\OO/\mm^L)-\rank M_\mu(s),\\
 \tau_k(f_s)&=\dim_\C(\OO/\mm^L)-\rank M_\tau(s).
 \end{aligned}
\]

Every matrix entry is polynomial in $s$.
For each matrix, choose a minor whose size equals its rank at
$s=0$ and whose value there is nonzero. If that rank is zero, no
condition is needed. Each chosen minor is polynomial in $s$ and
nonzero at $s=0$, so its nonvanishing defines a nonempty open
subset of the parameter line.

Intersect these sets with the ordinary locus. Each is nonempty and
open in the irreducible parameter line, so the intersection is
nonempty. On this intersection, both matrix ranks are at least
their values at $s=0$. Since the ambient dimension is fixed, the
corresponding quotient lengths are at most their values at $s=0$.
This proves
\eqref{eq:usc}.
\end{proof}

\subsection{Completion of the minimum theorem}

\begin{proof}[Proof of Theorems~\ref{thm:main} and~\ref{thm:family}]
Fix $f\in\Pall$ and $k\geq0$, and set $d=\deg f$.
Lemma~\ref{lem:usc} supplies an ordinary semihomogeneous member
$f_s\in\Pclass$ with $\mu_k(f)\geq\mu_k(f_s)$ and
$\tau_k(f)\geq\tau_k(f_s)$.
Proposition~\ref{prop:ordinarymu} gives $\mu_k(f_s)=B_{N,m}(k)$, and
Proposition~\ref{prop:semihomogeneous-values-bounds} gives the asserted
Tjurina lower bound for $f_s$, hence for $f$.
The parameter $s$ may depend on $k$; the resulting inequalities hold
for the fixed polynomial $f$ at every order.
\par
Every member of the family \eqref{eq:family} belongs to $\Pall$ and
has $\mu_k=B_{N,m}(k)$ by Proposition~\ref{prop:ordinarymu}. By
Corollary~\ref{cor:semihomogeneous-sharp}, the same polynomial attains
every Tjurina bound. Thus the lower bounds over $\Pall$ are the minima,
attained simultaneously for all $k\geq0$.
\end{proof}

\subsection{One open set for all orders}

\begin{proof}[Proof of Theorem~\ref{thm:generic}]
Let $U\subset V_{m,d}$ be the open locus
where $f_d\ne0$ and $f_m$ has isolated gradient.
It is nonempty and contained in $\Pclass$.
All its members have $\mu_k=B_{N,m}(k)$ for every $k$.

For each $f\in U$, put $A=\OO/J(f)$ and $a=[f]$.
Proposition~\ref{prop:ordinarymu} gives the uniform bound
$\mm_A^{S+1}=0$. If

$k\geq\max\{1,S-m+1\}$, then
\[
 a\mm_A^{k-1}\subseteq\mm_A^{m+k}=0
 \qquad\text{for every member of }U.
\]
The inequality $m+k>S$ also makes every summand of
$Q_{N,m}(k-1)$
zero. Thus the actual multiplication image and the correction in
the minimum formula vanish together. Equation \eqref{eq:defectk}
then gives the minimum automatically for these $k$, leaving only
$0\leq k\leq\max\{0,S-m\}$ to consider.

Fix one of the remaining orders $k$.

\emph{Step 1. $U_k$ is open.}
The uniform inclusion
$\mm^{S+1+k}\subseteq\mm^kJ(f)$ permits computation in
$\OO/\mm^{S+1+k}$. As in Lemma~\ref{lem:usc}, the image of
$(f,\mm^kJ(f))$ in this finite quotient is spanned by the columns of
a matrix $M_k(f)$ with entries polynomial in the coefficients of $f$.
Hence

\[
 \tau_k(f)=\dim_\C\bigl(\OO/\mm^{S+1+k}\bigr)-\rank M_k(f).
\]
Put $r_k=\max_{g\in U}\rank M_k(g)$ and define
\[
 U_k=\{f\in U:\rank M_k(f)=r_k\}.
\]
The complementary locus is cut out by the $r_k$-minors, so
$U_k$ is Zariski open in $U$. By construction, $  f\in U_k $ if and only if  
$\tau_k(f)=\min_{g\in U}\tau_k(g)$.

\emph{Step 2. The minimum on $U$ equals the full-class minimum.}
Take $\lambda=1$ in \eqref{eq:family}, with $g=0$ for $d=m+1$ and
$g=x_1^d$ for $d\geq m+2$. This polynomial belongs to $U$ and, by
Propositions~\ref{prop:ordinarymu} and~\ref{prop:attainment},
attains the minimum values of Theorem~\ref{thm:main}.

The polynomial just constructed belongs to every $U_k$.
Their finite intersection is therefore open in $V_{m,d}$ and

nonempty. On this intersection, the finitely many remaining Tjurina
numbers attain their minima; the larger orders have already been
covered by the nilpotence bound.

Every $\mu_k$ attains its minimum on all of $U$ by
Proposition~\ref{prop:ordinarymu}. This proves the theorem.
\end{proof}

\begin{remark}\label{rem:homogeneous}
If $d=m$, every member of $\Pclass$ is homogeneous and satisfies,
for $k\geq0$,
\begin{equation}\label{eq:homogeneous}
\mu_k(f)=B_{N,m}(k),\qquad \tau_k(f)=B_{N,m}(k)-\tbin{N+k-2}{N}.
\end{equation}
Indeed, Euler's identity gives $f\in J(f)$, so $[f]=0$ in the Milnor algebra.
The homogeneous Tjurina formula appears in
\cite[Corollary~4]{HYZHomogeneous}, and the corresponding
weighted homogeneous difference formula is given in
\cite[Proposition~2.6]{HLYZ}. In the present setting,
\eqref{eq:homogeneous} follows from Propositions~\ref{prop:ordinarymu}
and~\ref{prop:defect}.

\end{remark}

%%%%%%%%%%%%%%%%%%%%%%%%%%%%%%%%%%%%%%%
% examples
%%%%%%%%%%%%%%%%%%%%%%%%%%%%%%%%%%%%%%%
\section{Finite formulas and generating series of the minima}
\label{sec:explicit}
 
\subsection{Finite formulas}
 
The cumulative Hilbert function gives finite binomial expressions
for both $B_{N,m}$ and $Q_{N,m}$. Define

\begin{equation}\label{eq:H}
 H(v)=\sum_{i=0}^{v}h_i
     =\sum_{a=0}^{N}(-1)^a\binom Na\tbin{N+v-aq}{N}.
\end{equation}
We set $H(v)=0$ for $v<0$. Passing to the partial sums
amounts to dividing the series in \eqref{eq:ordinaryA} by $1-t$:
\[
 \sum_{v\geq0}H(v)t^v=\frac{\Hilb(\gr A,t)}{1-t}
       =\frac{(1-t^q)^N}{(1-t)^{N+1}}.
\]
The second equality in \eqref{eq:H} follows by expanding
the right-hand side: the numerator contributes
$(-1)^a\binom Na t^{aq}$,
and the coefficient of $t^{v-aq}$ in $(1-t)^{-(N+1)}$
is $\tbin{N+v-aq}{N}$. Thus $H(v)$ can be evaluated using
$N+1$ binomial terms.

The Milnor value $B_{N,m}(k)$ in \eqref{eq:B} also has the binomial
expression \cite[equation~(4.3)]{LiuYauZhang2026}

\begin{equation}\label{eq:Balternating}
 B_{N,m}(k)=q^N+\sum_{a=1}^{N}(-1)^{a+1}\binom Na
                    \tbin{N+k-1-(a-1)q}{N}.
\end{equation}

The splitting degree used in Lemma~\ref{lem:rank} also gives
closed formulas for the correction.
\begin{proposition}\label{prop:Qclosed}
Let $H$ be the cumulative Hilbert function \eqref{eq:H}, and put
\[
 c_j=\max\left\{j,\left\lceil\frac{S-m-1}{2}\right\rceil\right\}.
\]
Then:
\begin{enumerate}
\item for every $j\geq0$,
\begin{equation}\label{eq:Qclosed}
 Q_{N,m}(j)=H(c_j-1)-H(j-1)+H(S-c_j-m-1);
\end{equation}
\item at order zero,
\begin{equation}\label{eq:Qzero}
 Q_{N,m}(0)=
 H\left(\left\lceil\frac{S-m-1}{2}\right\rceil-1\right)
 +H\left(\left\lfloor\frac{S-m-1}{2}\right\rfloor\right).
\end{equation}
\end{enumerate}
\end{proposition}
\begin{proof}
\emph{(i).} Below $c_j$ the smaller Hilbert-function value is
$h_i$; from $c_j$ onward it is $h_{i+m+1}$. Hence
\[
 Q_{N,m}(j)=\sum_{i=j}^{c_j-1}h_i+\sum_{i=c_j+m+1}^{S}h_i.
\]
Symmetry, $h_i=h_{S-i}$, turns the second sum into an initial
segment: replacing $i$ by $S-v$ gives
\[
 \sum_{i=c_j+m+1}^{S}h_i
       =\sum_{v=0}^{S-c_j-m-1}h_v=H(S-c_j-m-1).
\]
The first sum is $H(c_j-1)-H(j-1)$, proving \eqref{eq:Qclosed}.

\emph{(ii).} If $S<m+1$, then $Q_{N,m}(0)=0$, and both
arguments of $H$ in \eqref{eq:Qzero} are negative. If $S\geq m+1$,
substituting $j=0$ in \eqref{eq:Qclosed} gives \eqref{eq:Qzero}.

\end{proof}

\subsection{Generating series of the minima}
Write $\mu_k^{\min}$ and $\tau_k^{\min}$ for the minima in
Theorem~\ref{thm:main}.  To describe all orders together, consider the generating functions
\[
 \mathcal M_{\min}(t)=\sum_{k\geq0}\mu_k^{\min}t^k,
 \qquad
 \mathcal T_{\min}(t)=\sum_{k\geq0}\tau_k^{\min}t^k.
\]
Put
\[
 P_{N,m}(t)=\sum_{a=1}^N(-1)^{a+1}\binom Na t^{(a-1)q}.
\]
By \cite[equation~(4.5)]{LiuYauZhang2026},
\begin{equation}\label{eq:Mseries}
 \mathcal M_{\min}(t)=\frac{q^N}{1-t}
          +\frac{tP_{N,m}(t)}{(1-t)^{N+1}}.
\end{equation}
The binomial term in the Tjurina defect is summed using the identity

\[
 \sum_{k\geq0}\tbin{N+k-b}{N}t^k
       =\frac{t^b}{(1-t)^{N+1}}\qquad(b\in\mathbb Z_{\geq0}).
\]

To sum the remaining correction $Q_{N,m}$, put

\[
 W_{N,m}(t)=\sum_{i\geq0}\min\{h_i,h_{i+m+1}\}t^i.
\]
Each rank in degree $i$ contributes to $Q_{N,m}(0),\ldots,Q_{N,m}(i)$.
Summing these contributions in the other order gives
\begin{align*}
 \sum_{j\geq0}Q_{N,m}(j)t^j
   &=\sum_{i\geq0}\min\{h_i,h_{i+m+1}\}(1+t+\cdots+t^i)\\
   &=\frac{Q_{N,m}(0)-tW_{N,m}(t)}{1-t}.
\end{align*}
The minimum formula uses $Q_{N,m}(0)$ at order zero and $Q_{N,m}(k-1)$
at orders $k\geq1$. These contributions give the total correction
series
\[
 Q_{N,m}(0)+t\sum_{j\geq0}Q_{N,m}(j)t^j
       =\frac{Q_{N,m}(0)-t^2W_{N,m}(t)}{1-t}.
\]
The homogeneous correction $\tbin{N+k-2}{N}$ contributes
$t^2/(1-t)^{N+1}$ by the same shifted identity with $b=2$.
Subtracting both contributions from \eqref{eq:Mseries} gives
\begin{equation}\label{eq:Tseries}
 \mathcal T_{\min}(t)=\mathcal M_{\min}(t)
     -\frac{t^2}{(1-t)^{N+1}}
     -\frac{Q_{N,m}(0)-t^2W_{N,m}(t)}{1-t}.
\end{equation}
At each prescribed degree $d>m$, the minimum Milnor and Tjurina
series are again $\mathcal M_{\min}(t)$ and $\mathcal T_{\min}(t)$.
At $d=m$, the Milnor series is unchanged, while the Tjurina series is
$\mathcal M_{\min}(t)-t^2/(1-t)^{N+1}$ by
Remark~\ref{rem:homogeneous}.

Since $W_{N,m}(1)=Q_{N,m}(0)$, the last fraction in
\eqref{eq:Tseries} is a polynomial, in agreement with
the eventual vanishing of $Q_{N,m}$ noted after \eqref{eq:Q}.

\section{Boundary cases and examples}\label{sec:examples}

\subsection{Boundary cases}

\emph{One variable.}
For $N=1$, write $f=x^m u(x)$ with $u(0)\ne0$.
Then $J(f)=(x^{m-1})$, and every member of the class has
\[
 \mu_k=m+k-1,\qquad \tau_0=m-1,\qquad \tau_k=m\quad(k\geq1).
\]
This agrees with Theorem~\ref{thm:main} because $Q_{1,m}\equiv0$.

\emph{Multiplicity two.}
For $m=2$ in any dimension, $h_0=1$ and all other $h_i$
vanish, so the minima are
\[
 \mu_k^{\min}=\binom{N+k}{N},\qquad
 \tau_k^{\min}=\binom{N+k}{N}-\tbin{N+k-2}{N}.
\]
Here the nonhomogeneous correction vanishes for every $k$ and every
allowed degree $d$.

\emph{The first higher order.}
For every $N,m$, \eqref{eq:Balternating} gives
$B_{N,m}(1)=q^N+N$, and consequently
\begin{equation}\label{eq:orderone}
 \tau_1^{\min}=\tau_0^{\min}+N.
\end{equation}
The same correction $Q_{N,m}(0)$ occurs at both orders; the difference comes
entirely from $B_{N,m}(1)-B_{N,m}(0)=N$.

\subsection{Plane curve singularities}

For $N=2$, the Hilbert function begins with $h_i=i+1$ for
$0\leq i<q$.

For $m\geq5$, we have $S-m-1=m-5$, so symmetry gives
\[
 Q_{2,m}(0)=\sum_{i=0}^{m-5}\min\{i+1,m-4-i\}
     =\left\lfloor\frac{(m-3)^2}{4}\right\rfloor.
\]
To see the last equality, pair terms from the two ends of the sum.
If $m-3=2a$, the sum is $2(1+\cdots+(a-1))+a=a^2$.
If $m-3=2a+1$, it is $2(1+\cdots+a)=a(a+1)$.
For $m=2,3,4$, we have $S<m+1$, so $Q_{2,m}(0)=0$,
again agreeing with this expression. Hence, over $\mathcal P_{2;m}$,
\begin{equation}\label{eq:plane}
 \min\tau_0=(m-1)^2-\left\lfloor\frac{(m-3)^2}{4}\right\rfloor
 \qquad(N=2).
\end{equation}
This is equivalently $\lfloor(3m^2-2m-4)/4\rfloor$, the
classical minimum for ordinary plane multiple points
\cite{BGM,CaninoGimiglianoIda}.
Canino--Gimigliano--Id\`a
\cite[Theorem~5.7 of arXiv v2]{CaninoGimiglianoIda} give explicit attainers:
$x^m+y^m+x^{m/2+1}y^{m/2}$ for even $m$, and
$x^m+y^m+x^{(m+1)/2}y^{(m+1)/2}$ for odd $m$.

Eriksson \cite{ErikssonFloors} also describes an Euler-reduced
generator count giving the same lower bound for arbitrary isolated
plane germs.

\subsection{Surface singularities}

In three variables, $S-m-1=2m-7$, and both terms in
\eqref{eq:Qzero} are $H(m-4)$. When $m\geq4$, all degrees through
$m-4$ are below $q$, so the monomial complete intersection has
the same components as the polynomial ring in this range. Thus
\[
 H(m-4)=\sum_{v=0}^{m-4}\binom{v+2}{2}=\binom{m-1}{3}.
\]
For $m=2,3$, $H(m-4)$ and the corresponding truncated binomial
both vanish. Consequently $Q_{3,m}(0)=2\tbin{m-1}{3}$ for every $m\geq2$,
and
\begin{equation}\label{eq:surface}
 \min\tau_0=(m-1)^3-2\tbin{m-1}{3}\qquad(N=3).
\end{equation}
For $m\geq4$, this recovers the values of Wahl and Eriksson
\cite[Example~4.7]{Wahl},
\cite[Theorem~2.2 and Proposition~3.2]{Eriksson} discussed in the
introduction.

Table~\ref{tab:examples} lists the first five orders in dimensions two,
three, and four.

\begin{table}[ht]
\centering
\renewcommand{\arraystretch}{1.15}
\caption{Sharp minima over $\Pall$. Each sequence lists the
values at $k=0,1,2,3,4$; they also apply at any fixed degree $d>m$.}
\label{tab:examples}
\begin{tabular}{ccc}
\toprule
$(N,m)$ & $(\mu_k^{\min})_{k=0}^4$ & $(\tau_k^{\min})_{k=0}^4$\\
\midrule
$(2,5)$ & $16,18,22,28,36$ & $15,17,21,25,30$\\
$(3,4)$ & $27,30,39,57,84$ & $25,28,37,53,74$\\
$(4,3)$ & $16,20,36,70,126$ & $15,19,35,65,111$\\
$(4,4)$ & $81,85,101,141,215$ & $71,75,91,131,199$\\
\bottomrule
\end{tabular}
\end{table}
\FloatBarrier

\subsection{A worked example}
Write $\ell=x_1+\cdots+x_N$. Taking $\lambda=1$ and $g=0$ in
Theorem~\ref{thm:family} gives the simultaneous attainer

\begin{equation}\label{eq:unrestrictedattainer}
 F_{N;m}=\sum_{i=1}^N x_i^m+\ell^{m+1}.
\end{equation}
Consider $f=x^4+y^4+z^4+(x+y+z)^5$, the specialization of $F_{N;m}$
at $(N,m)=(3,4)$. By Theorem~\ref{thm:family}, it attains every minimum
over $\mathcal P_{3;4}$, and adding terms of order at least six preserves
these values. Both correction terms in \eqref{eq:minnuk} are nonzero at
$k=2$. The graded Milnor algebra of $f$ has Hilbert function

\[
 (h_0,\ldots,h_6)=(1,3,6,7,6,3,1).
\]
As computed in Example~\ref{ex:rank}, $Q_{3,4}(0)=2$,
$Q_{3,4}(1)=1$, and $Q_{3,4}(j)=0$ for $j\geq2$.
Over $\mathcal P_{3;4}$, the first three minimum values are
\[
 (\mu_0^{\min},\mu_1^{\min},\mu_2^{\min})=(27,30,39),\qquad
 (\tau_0^{\min},\tau_1^{\min},\tau_2^{\min})=(25,28,37).
\]
At $k=2$, the Milnor formula reads
\[
 B_{3,4}(2)=\binom{7}{3}+\sum_{i=5}^{6}h_i,
\]
where $\binom{7}{3}=35$ counts the monomials of degree below five
and $h_5+h_6=3+1=4$ is the dimension of the remaining graded
Milnor components. Hence

\[
 B_{3,4}(2)=35+4=39.
\]
The homogeneous Tjurina correction is $\binom{3}{3}=1$, giving
$39-1=38$. The additional correction $Q_{3,4}(1)=1$ lowers this value to
$37$. Thus \eqref{eq:minnuk} separates the decrease already present
in the homogeneous case from the further decrease produced by the
perturbation.

For this example, the perturbation lowers the Tjurina number only at
orders $k=0,1,2$. From $k=3$ onward, the value is
$B_{3,4}(k)-\tbin{k+1}{3}$, the same as in the homogeneous case.
In general, for $k\geq1$, the nonhomogeneous correction vanishes as
soon as $m+k>S$. 

\bibliographystyle{amsplain}
\bibliography{refs}

\end{document}